\documentclass[11pt]{article}
\usepackage[a4paper, margin=2.5cm]{geometry}
\usepackage{amsmath,amsthm,amsfonts,amssymb}
\usepackage{mathpazo}
\usepackage[authoryear]{natbib}
\usepackage{hyperref}
\usepackage{booktabs}
\usepackage{multirow}
\usepackage{tabularx}
\usepackage{url}
\usepackage[skip=.5\baselineskip plus 2pt, indent=20pt]{parskip}
\usepackage{setspace}
\newtheoremstyle{thmstyle}{\topsep}{\topsep}{\itshape}{}{\bfseries}{.}{.5em}{}
\theoremstyle{thmstyle}
\newtheorem{theorem}{Theorem}[section]
\newtheorem{proposition}{Proposition}[section]
\newtheorem{corollary}{Corollary}[section]
\newtheorem{definition}{Definition}[section]

\usepackage{mathtools}

\begin{document}

\title{Toll Allocation in General Transportation Networks}
\author{Alejandro Saavedra-Nieves\\[4pt]
\small CITMAga, MODESTYA Research Group\\
\small Departamento de Estat\'istica, An\'alise Matem\'atica e Optimizaci\'on\\
\small Universidade de Santiago de Compostela, Spain\\
\small \href{mailto:alejandro.saavedra.nieves@usc.es}{alejandro.saavedra.nieves@usc.es}}
\date{}
\maketitle

\begin{abstract}
The allocation of transportation revenues among the infrastructure
segments used by travellers is a relevant problem for the economic
management of transportation networks.  Existing approaches, such as
the Segments Equal Sharing (SES) rule, address this problem in linear
highway structures, but do not explicitly account for the more complex
topologies found in general networks.  This paper extends the highway
toll allocation problem to general networks and introduces the
Generalized Segments Equal Sharing (GSES) rule.  Under this rule, the
revenue generated by each journey is distributed equally among the
segments belonging to the path used, and the allocation to each
segment is obtained by aggregating its shares across all
journeys.  We provide two axiomatic characterizations of the proposed
rule based on natural properties in this context.  We also establish
that GSES coincides with the Shapley value of an associated segments
allocation game, providing a cooperative game-theoretic foundation
for the rule.
\end{abstract}


\section{Introduction}\label{sec:introduction}

The economic management of transportation infrastructure requires not
only determining the revenues generated by its use, but also
understanding how those revenues can be attributed to the individual
components of the network.  In tolled systems, users pay for journeys
that may traverse several infrastructure segments; the resulting
revenues must therefore be distributed among those segments.  A
consistent allocation can support maintenance planning, investment
prioritization, and performance evaluation~\citep{ECMT2002,
dePalma2007, Weinreich2021}.

The allocation problem becomes non-trivial in general networks, where
branches, intersections, and alternative paths allow several routes to
connect the same origin--destination pair~\citep{Sheffi1985,
OrtuzarWillumsen2011}.  Cooperative game theory provides a natural
framework for this class of problems: by treating infrastructure
segments as players, one can relate the revenues generated by journeys
to the contributions of the segments involved.  This approach has been
applied to airport landing fees~\citep{Littlechild1977, Vazquez1997,
SaavedraNieves2026}, infrastructure cost allocation in line-based
systems~\citep{Estan2021, Babaei2022}, revenue sharing among
transport providers~\citep{Ding2024}, and highway cost
allocation~\citep{Dong2012, Kuipers2013, Sudholter2017, Gomez2024}.

Most directly related to our work, \citet{Wu2024} introduced the
highway toll allocation problem for one-way linear highways and
proposed the Segments Equal Sharing (SES) rule, under which the toll
paid by a user is distributed equally among the segments traversed.
They showed that SES coincides with the Shapley value~\citep{Shapley1953}
of an associated TU-game and provided axiomatic characterizations.
\citet{Soto2026} subsequently studied alternative allocation rules in
the same linear setting.

The linear structure underlying these contributions restricts the
class of network topologies that can be represented.  In particular,
it does not address the allocation problem arising in networks with
branches, intersections, shared segments, and alternative paths.  To
the best of our knowledge, the present paper provides the first
extension of the highway toll allocation problem and its equal-sharing
principle to general network topologies.

Our contribution is threefold.  First, we formulate the toll
allocation problem for general networks and introduce the Generalized
Segments Equal Sharing (GSES) rule as a path-based generalization of
SES.  Second, we provide two axiomatic characterizations of GSES.
Third, we prove that GSES coincides with the Shapley value of an
associated segments allocation game and that this allocation always
belongs to the core, connecting the rule to well-established
cooperative solution concepts.

The paper is organized as follows.  Section~\ref{sec:prel} reviews
the required cooperative game theory concepts and the linear highway
model.  Section~\ref{sec:generalproblem} formulates the general
network model.  Section~\ref{sec:GESrule} introduces GSES, develops
its axiomatic characterizations, and establishes the connection with
the Shapley value and the core.  Section~\ref{sec:concluding}
concludes.

\section{Preliminaries}\label{sec:prel}

\subsection{Cooperative game theory}

A \emph{TU-game} is a pair $(N,v)$ where $N$ is a finite player set
and $v:2^N\to\mathbb{R}$ is a characteristic function with
$v(\emptyset)=0$.  The \emph{core}~\citep{Gillies1953} is
\[
C(N,v)=\bigl\{x\in\mathbb{R}^N:\textstyle\sum_{i\in N}x_i=v(N),\;
\sum_{i\in S}x_i\geq v(S)\text{ for every }S\subsetneq N\bigr\}.
\]
The \emph{Shapley value}~\citep{Shapley1953} assigns each player its
expected marginal contribution over uniformly random orderings:
\[
\operatorname{Sh}_i(N,v)=
\sum_{S\subseteq N\setminus\{i\}}
\frac{|S|!\,(|N|-|S|-1)!}{|N|!}
\bigl[v(S\cup\{i\})-v(S)\bigr].
\]

\subsection{The linear highway model and SES rule}

\citet{Wu2024} consider a one-way linear highway with segment set
$N=\{1,\ldots,n\}$ ordered consecutively.  For any $h\leq k$, the
trip $[h,k]$ traverses segments $\mathcal{P}([h,k])=\{h,h+1,\ldots,k\}$
and generates toll revenue $t_{hk}\geq 0$.  A \emph{toll matrix}
$T=(t_{hk})\in\mathbb{R}_+^{n\times n}$ is upper triangular (with
$t_{hk}=0$ for $h>k$).

\begin{definition}[\citealt{Wu2024}]
The \emph{Segments Equal Sharing} (SES) method is defined, for every
$T\in\mathcal{T}^N$ and every $i\in N$, by
\[
f^{\mathrm{Se}}_i(N,T)=
\sum_{h=1}^{i}\sum_{k=i}^{n}\frac{t_{hk}}{k-h+1}.
\]
\end{definition}

The \emph{segments allocation game} $(N,v^{\mathrm{Se}})$ assigns
to each coalition $S\subseteq N$ the revenue generated by trips whose
segments all belong to $S$:
$v^{\mathrm{Se}}(S)=\sum_{h,k:\,\mathcal{P}([h,k])\subseteq S}t_{hk}$.
\citet{Wu2024} showed that $f^{\mathrm{Se}}(N,T)=\operatorname{Sh}(N,v^{\mathrm{Se}})$
and provided two axiomatic characterizations of SES.

\section{Toll Allocation in General Networks}\label{sec:generalproblem}

Let $N$ be a set of $n$ segments.  For any pair of nodes $h$ and $k$,
let $\mathcal{P}([h,k])$ denote the set of admissible paths from $h$
to $k$.  Each path $p_\ell$ is characterized by its set of segments
$\mathcal{P}(p_\ell)\subseteq N$.  Every path contains at least one
segment; a single segment is a path.  The toll collected from users
travelling from $h$ to $k$ along $p_\ell$ is $t_{hk}^\ell\geq 0$.  A
\emph{toll profile} is the collection
$T=(t_{hk}^\ell)_{h,k,\,p_\ell\in\mathcal{P}([h,k])}$; the set of all
such profiles is $\mathcal{T}^N$.  The total revenue is
$\tau(T)=\sum_{h,k}\sum_{p_\ell\in\mathcal{P}([h,k])}t_{hk}^\ell$.

A \emph{highway toll allocation problem on a general network} is a
tuple $(N,T)$.  An \emph{allocation} is a vector
$x\in\mathbb{R}_+^N$, and a \emph{toll allocation method} is a
mapping $f$ assigning an allocation to each such problem.

When the network is a linear highway with a single path between every
pair of nodes, the toll profile reduces to a toll matrix and the
general framework specializes to that of \citet{Wu2024}.

\begin{definition}
For every $T\in\mathcal{T}^N$, the \emph{segments allocation game}
$(N,v^{\mathrm{GSe}})$ is defined by
\[
v^{\mathrm{GSe}}(S)=
\sum_{h,k\in N}
\sum_{\substack{p_\ell\in\mathcal{P}([h,k])\\\mathcal{P}(p_\ell)\subseteq S}}
t_{hk}^\ell,
\qquad S\subseteq N.
\]
\end{definition}

\section{The GSES Method}\label{sec:GESrule}

\subsection{Definition}

\begin{definition}\label{def:GSES}
For every toll allocation problem $(N,T)$, the
\emph{Generalized Segments Equal Sharing} (GSES) method is defined by
\[
f_i^{\mathrm{GSES}}(N,T)=
\sum_{h,k\in N}
\sum_{\substack{p_\ell\in\mathcal{P}([h,k])\\i\in\mathcal{P}(p_\ell)}}
\frac{t_{hk}^\ell}{|\mathcal{P}(p_\ell)|},
\qquad i\in N.
\]
\end{definition}

Thus, the revenue associated with each path is shared equally among
the segments composing that path, and each segment's total allocation
aggregates its shares across all revenue-generating paths.  In the
linear setting, $|\mathcal{P}(p_\ell)|=k-h+1$ for the unique path
from $h$ to $k$, and GSES reduces to SES.

\subsection{Axiomatic characterizations}

We introduce the following properties for a toll allocation method $f$.

\begin{itemize}
\item \textbf{Efficiency (EFF).}
$\sum_{i\in N}f_i(N,T)=\tau(T)$ for every $T\in\mathcal{T}^N$.

\item \textbf{Additivity (ADD).}
$f(N,T+T')=f(N,T)+f(N,T')$ for every $T,T'\in\mathcal{T}^N$.

\item \textbf{Toll linearity (TL).}
$f(N,bT+b'T')=bf(N,T)+b'f(N,T')$ for every $T,T'\in\mathcal{T}^N$
and $b,b'\geq 0$.

\item \textbf{Inessential segment (IS).}
If $t_{hk}^\ell=0$ for every path $p_\ell$ containing $i$, then
$f_i(N,T)=0$.

\item \textbf{Equal sharing within a path (ES).}
For every $[h,k]$ and $p_\ell\in\mathcal{P}([h,k])$,
\[
f_i(N,\Delta^{hk,\ell})=
\begin{cases}
1/|\mathcal{P}(p_\ell)|, & i\in\mathcal{P}(p_\ell),\\
0, & \text{otherwise},
\end{cases}
\]
where $\Delta^{hk,\ell}$ assigns one unit of revenue to $p_\ell$ and
zero to all other paths.

\item \textbf{Necessary segment (NS).}
If every path with positive toll revenue contains $i$, then
$f_i(N,T)\geq f_j(N,T)$ for every $j\in N$.
\end{itemize}

Property (TL) is introduced in~\citet{Soto2026} for the linear
setting.  Properties (ADD), (EFF), and (IS) extend the corresponding
axioms of~\citet{Wu2024} to general networks.  Property (ES) extends
their segment symmetry axiom to the path level.  Property (NS) extends
their necessary-segment axiom.

\begin{theorem}\label{th:1ca}
The GSES method is the unique toll allocation method satisfying
toll linearity \emph{(TL)} and equal sharing within a path \emph{(ES)}.
\end{theorem}

The proof proceeds in two steps.  Existence follows from the
definition: GSES satisfies (TL) since each $t_{hk}^\ell$ enters
linearly, and satisfies (ES) by direct substitution.  Uniqueness
follows because (TL) allows any toll profile to be decomposed into
unit profiles $\Delta^{hk,\ell}$, and (ES) then uniquely pins down
the allocation on each unit profile.

\begin{theorem}\label{th:2c}
The GSES method is the unique toll allocation method satisfying
additivity \emph{(ADD)}, efficiency \emph{(EFF)}, the inessential
segment property \emph{(IS)}, and the necessary segment property
\emph{(NS)}.
\end{theorem}

For uniqueness, consider any single-path profile $t\Delta^{hk,\ell}$.
Since every segment of $p_\ell$ is necessary, (NS) forces all of them
to receive equal allocations; (IS) assigns zero to segments outside
$p_\ell$; (EFF) then yields $t/|\mathcal{P}(p_\ell)|$ each.
Extending to arbitrary profiles by (ADD) recovers the GSES formula.

The logical independence of the axioms in each characterization is
verified in the online appendix through explicit counterexamples.

\section{The Shapley Value and the Core}\label{sec:shapley}

\begin{theorem}\label{thm:gses-shapley}
Let $(N,T)$ be a highway toll allocation problem on a general network.
Then
\[
f^{\mathrm{GSES}}(N,T)=\operatorname{Sh}(N,v^{\mathrm{GSe}}).
\]
\end{theorem}

The characteristic function $v^{\mathrm{GSe}}$ decomposes as a
non-negative linear combination of unanimity games,
$v^{\mathrm{GSe}}=\sum_{h,k}\sum_{p_\ell}t_{hk}^\ell\,
u_{\mathcal{P}(p_\ell)}$.  The additivity of the Shapley value and
the fact that $\operatorname{Sh}_i(N,u_S)=1/|S|$ for $i\in S$ and
$0$ otherwise then give
$\operatorname{Sh}_i(N,v^{\mathrm{GSe}})=
\sum_{h,k}\sum_{i\in\mathcal{P}(p_\ell)}t_{hk}^\ell/|\mathcal{P}(p_\ell)|
=f_i^{\mathrm{GSES}}(N,T)$.

\begin{proposition}\label{prop:convex}
For every $T\in\mathcal{T}^N$, the game $(N,v^{\mathrm{GSe}})$ is
convex.
\end{proposition}

Convexity follows because every unanimity game is convex and
convexity is preserved under non-negative linear combinations.  Since
the Shapley value of a convex game belongs to the core
\citep{shapley1971cores}, Theorem~\ref{thm:gses-shapley} and
Proposition~\ref{prop:convex} immediately yield:

\begin{corollary}\label{cor:core}
For every $(N,T)$, $f^{\mathrm{GSES}}(N,T)\in C(N,v^{\mathrm{GSe}})$.
\end{corollary}

This extends the core membership result of~\citet{Wu2024} beyond
linear highways without any topological assumption on the network.

As in~\citet{Soto2026}, GSES belongs to a broader family of
allocation methods
$\mathcal{F}(N,T)=\{f^\alpha(N,T):\alpha\in\mathcal{A}(N,T)\}$
where $f_i^\alpha(N,T)=\sum_{h,k}\sum_{i\in\mathcal{P}(p_\ell)}
\alpha_{hk}^{i,\ell}t_{hk}^\ell$ and the weights
$\alpha_{hk}^{i,\ell}\geq 0$ are arbitrary.  The GSES weights are
$\alpha_{hk}^{i,\ell}=1/|\mathcal{P}(p_\ell)|$.

\begin{proposition}\label{prop:core-family}
The core $C(N,v^{\mathrm{GSe}})$ equals the set of allocations in
$\mathcal{F}(N,T)$ whose weights satisfy
$\sum_{i\in\mathcal{P}(p_\ell)}\alpha_{hk}^{i,\ell}=1$ for every
$h,k$ and every $p_\ell\in\mathcal{P}([h,k])$.
\end{proposition}

The proof uses the convexity of $(N,v^{\mathrm{GSe}})$ and an
extreme-point argument ordering segments by an arbitrary permutation,
which does not require the segments of a path to form a contiguous
interval.  Corollary~\ref{cor:core} is recovered as a special case
since the GSES weights trivially sum to one on every path.

\section{Concluding Remarks}\label{sec:concluding}

We have extended the highway toll allocation problem of~\citet{Wu2024}
from linear highways to general transportation networks.  The proposed
GSES method distributes the revenue generated by each path equally
among its constituent segments.  It is characterized axiomatically in
two ways: through toll linearity and equal sharing within a path, and
through additivity, efficiency, the inessential segment property, and
the necessary segment property.  It also coincides with the Shapley
value of the associated segments allocation game and always belongs to
its core.

The framework applies whenever usage data can be organized into
origin--destination trips and their underlying routes.  An empirical
illustration using passenger validation data from the Hangzhou metro
network (2019) is available in the full version of the paper.

Several directions remain open.  First, the sensitivity of the
allocation to the route-reconstruction procedure in networks where
multiple shortest paths exist deserves further study.  Second,
identifying analogues of the other allocation methods studied
in~\citet{Soto2026} for general networks is non-trivial, as those
rules rely on the linear ordering of segments along a trip.  Third,
applications to other transportation contexts and datasets would help
establish the practical scope of the method.

\section*{Acknowledgements}
This work is part of project PID2021-124030NB-C32, funded by
MCIN/AEI/10.13039/501100011033 and ``ERDF A way of making Europe''/EU.
Research also funded by Grupos de Referencia Competitiva ED431C~2021/24
from the Conseller\'ia de Cultura, Educaci\'on e Universidades, Xunta
de Galicia, Spain.

\bibliographystyle{plainnat}

\end{document}